\documentclass[11pt,reqno]{amsart}

\usepackage[foot]{amsaddr}
\usepackage{nicefrac, xfrac}
\usepackage{enumerate}
\usepackage{amsfonts, amsmath, amssymb, amscd, amsthm, bm}

\usepackage{enumerate}
\usepackage{url}
\usepackage[linktocpage=true,colorlinks,citecolor=magenta,linkcolor=blue,urlcolor=magenta]{hyperref}
\usepackage{multicol}
\usepackage{comment}
\usepackage{xcolor}
\usepackage{mathrsfs, graphicx, subcaption, caption}

 \newtheorem{theorem}{Theorem}[section]
 \newtheorem{corollary}[theorem]{Corollary}
 
 \newtheorem{proposition}[theorem]{Proposition}
 \theoremstyle{definition}
 \newtheorem{definition}[theorem]{Definition}
 \newtheorem{remark}[theorem]{Remark}
 \newtheorem*{ack}{Acknowledgments}
 
 \newtheorem{problem}[theorem]{Problem}
 \numberwithin{equation}{section}
\def\R{\mathbb R}

\def\p#1{\partial #1}

\numberwithin{equation}{section}

\newcounter{rom}
\renewcommand{\therom}{(\roman{rom})}

\makeatletter
\@namedef{subjclassname@2020}{\textup{2020} Mathematics Subject Classification}
\makeatother

\title[symmetry of hypersurfaces]{Symmetry of hypersurfaces with symmetric boundary}

\begin{document}

\author{Hui Ma$^{1}$}
\address{$^1$Department of Mathematical Sciences, Tsinghua University,
Beijing 100084, P.R. China} 
\email{ma-h@mail.tsinghua.edu.cn}

\author{Chao Qian$^2$}
\address{$^2$School of Mathematics and Statistics, Beijing Institute of Technology, Beijing 100081, P.R. China}
\email{6120150035@bit.edu.cn}
\author{Jing Wu$^1$}
\email{w3024v@gmail.com}

\author{Yongsheng Zhang$^3$}
\address{$^3$Academy for Multidisciplinary Studies, Capital Normal University, Beijing, 100048, P.R. China}
\email{yongsheng.chang@gmail.com}

\begin{abstract}
Let $G$ be a compact connected subgroup of $SO(n+1)$. In $\mathbb{R}^{n+1}$,
we gain interior $G$-symmetry for minimal hypersurfaces and hypersurfaces of constant mean curvature (CMC) which have  $G$-invariant boundaries and $G$-invariant contact angles along boundaries. The main ingredients of the proof are to build an associated Cauchy problem based on infinitesimal Lie group actions, and to apply Morrey’s regularity theory and the Cauchy-Kovalevskaya Theorem. Moreover, we also investigate the same kind of symmetry inheritance from boundaries for hypersurfaces of constant higher order mean curvature and Helfrich-type hypersurfaces in $\mathbb{R}^{n+1}$. 
\end{abstract}



\maketitle

\section{Introduction}\label{sec1}
The classical Plateau problem aims to find a minimal surface spanned by a given closed curve. 
Whether a solution can inherit the symmetry from its boundary has been studied in the area-minimizing integral currents setting, e.g., Lawson \cite{Lawson72}, Bindschadler \cite{Bind80}, Morgan \cite{Morgan86}, Lander \cite{Lander88} and etc. 
More generally, an intriguing problem is 
\begin{problem}\label{P1}
What kind of boundary symmetry can pass to the interior of certain interesting hypersurfaces in $\mathbb{R}^{n+1}$?
\end{problem}

A famous open conjecture states that a  compact, embedded surface with non-zero constant mean curvature and circular boundary must be a spherical cap (or equivalently a surface of revolution).
Under additional conditions, Koiso \cite{Koiso86} and Earp-Brito-Meeks-Rosenberg \cite{EBMR91} proved the conjecture. An analogous conjecture for hypersurfaces of constant higher order mean curvature was settled by Al\'{\i}as and Malacarne \cite{AM02}. For minimal hypersurfaces,  Schoen \cite{Schoen83} showed certain reflection symmetries can be induced from the boundary to the interior. In all these papers, the key idea is to apply the Alexandrov reflection method, which shows its power in the case with spherical boundary.

However, the Alexandrov reflection method seems unapplicable to symmetry given by general Lie subgroup $G\subset SO(n+1)$. 
Inspired by the work of Palmer and P\'ampano \cite{PP21} for CMC surfaces with circular boundary in $\mathbb{R}^3$, we start explorations on 
non-spherical boundary of higher dimension in this paper.
We study symmetries  
by a compact connected
Lie subgroup $G$ (not necessarily of cohomogeneity 2 in $\R^{n+1}$) of $SO(n+1)$ 
and those from non-homogeneous isoparametric foliations of $S^n(1)$.
We focus on minimal hypersurfaces, CMC hypersurfaces, hypersurfaces of constant higher order mean curvature and Helfrich-type hypersurfaces in $\mathbb{R}^{n+1}$.
The uniform idea is to translate the invariance problem with constant contact angle (see Definition \ref{defAngle}) along a boundary component to an elliptic Cauchy problem (see \eqref{Cauchy2}) and then combine Morrey's regularity theory with the Cauchy-Kovalevskaya Theorem for local $G$-invariance of hypersurfaces. Furthermore the invariance can extend to global $G$-invariance by assuming connectedness and completeness of  hypersurfaces. 
Generally, these two assumptions are indispensable to derive the global $G$-invariance.

Let us introduce the concept of contact angle  which we shall focus on.

\begin{definition}\label{defAngle}
Let $X: \Sigma \rightarrow \mathbb{R}^{n+1}$ be an immersion of a hypersurface with boundary $\partial \Sigma$. 
When a connected component $\Gamma$ of $\partial\Sigma$ lies in a hypersphere $S^n(R)$,
 one can define \emph{contact angle} $\theta \in [0,2\pi)$ 
 pointwise on $\Gamma$  by
\begin{equation}\label{normalDEF}
\left(\begin{array}{cc}
\nu   \\
 \mathfrak{n}  \\
 \end{array}\right)=
\left( \begin{array}{cc}
\cos \theta& \sin\theta  \\
-\sin\theta & \cos \theta \\
\end{array} \right)
\left(\begin{array}{cc}
\frac{X}{R}  \\
 N \\
 \end{array}\right),
 \end{equation}
where
$\nu$ is the unit normal of $\Sigma$ in $
\mathbb{R}^{n+1}$ defined along $\partial \Sigma$, $\mathfrak{n}$ is the outward unit conormal of $\partial\Sigma$ in $\Sigma$ and
$N$ the unit conormal of $\Gamma$ in $S^n(R)$.
\end{definition}

Throughout our paper, for simplicity
we assume that $X:\Sigma \to\mathbb{R}^{n+1}$ is an embedding or immersion of a connected oriented hypersurface with boundary $\partial\Sigma$. If there is no ambiguity in the embedding case, we will identify $X(\Sigma)$ with $\Sigma$ and $X(\partial\Sigma)$ with $\partial\Sigma$, respectively.

As one of our main results, we obtain
\begin{theorem}\label{gen2}
Let $G$ be a compact connected Lie subgroup  of $SO(n+1)$ and  $X:\Sigma \rightarrow \mathbb{R}^{n+1}$ an embedding of a connected hypersurface with boundary $\partial\Sigma$. Assume that there exists  a $G$-invariant connected component $\Gamma$ of $\partial \Sigma$ which is a real analytic submanifold lying in $S^n(R)$. 
If the contact angle $\theta$ along $\Gamma$ is $G$-invariant 
and $\Sigma$ is one of the following cases:
\begin{itemize}
\item[\rm{(i)}] a $C^1$ minimal hypersurface;
\item[\rm{(ii)}] a $C^1$ hypersurface of constant mean curvature;
\item[\rm{(iii)}] a $C^{2}$ hypersurface of constant $r$-th mean curvature that contains an interior elliptic point, satisfying $r>1$ and $H_r>0$ with respect to a unit normal vector field;
\item[\rm{(iv)}]  a $C^{4,\alpha}$ Helfrich-type\footnote{See Definition \ref{defHelfrich}.} hypersurface,
\end{itemize}
then the interior of $\Sigma$ is locally $G$-invariant\footnote{See Definition \ref{def2.3}.}.

If in addition $\Sigma$ is complete with respect to the induced metric and  $\partial\Sigma$ is $G$-invariant, then $\Sigma$ is $G$-invariant.
\end{theorem}

\begin{remark}\label{rk1}
$G$-invariance condition implies that $\Gamma$ is a union of $G$-orbits, along each of which the contact angle $\theta$ is constant. 
This actually is crucial.
    For instance, consider two circles of the same size in a pair of parallel planes in $\mathbb{R}^3$. When they are close enough and share the same rotation axis, the soap film formed between them is a catenoid. However, if we displace  one of the circles slightly inside the plane, then the  soap film loses its rotational symmetry, even locally. 
\end{remark}

\begin{remark}\label{rk2}
For $G$ of high cohomogeniety, the assumption that $\partial \Sigma$ is $G$-invariant does not imply that
each connected component of $\partial \Sigma$ is real analytic.
Namely, boundary components other than $\Gamma$ in the conclusion may not be real analytic.
\end{remark}

\begin{remark}\label{rk3}
Actually, 
each $\Sigma$ among $\rm{(i)}$-$\rm{(iv)}$ is a hypersurface 
with interior analyticity.
Take $G$-invariant boundary component $\Gamma$ for example. 
As our method is to apply
the Cauchy-Kowalevskaya Theorem,
if, along some open set $\underline{U}\subset\Gamma$ (not necessarily along the entire $\Gamma$), $\underline{U}$ is real analytic and both  $\underline{U}$ and $\theta$
are locally $G$-invariant,
then the same interior local $G$-invariance holds (see Theorem \ref{localthm} for details).
This will enable us to include many more boundary situations, cf. Remarks \ref{rk1} and \ref{rk2}.
\end{remark}

\begin{remark}
 For area-minimizing integral current of arbitrary codimension, Lander \cite{Lander88} showed the invariance of an area-minimizing integral current with boundary $B$ under a polar group action of $G$, where $B$ is supposed to be $G$-invariant and lying in the union of the principal orbits. Unlike Lander's requirement, in Theorem \ref{gen2} the orbits in $\Gamma$ are not necessarily of same type and the action of $G$ is not necessarily a polar action.
\end{remark}

When the contact angle equals to $\frac{\pi}{2}$, we have the following immediate application of Theorem \ref{gen2}.

\begin{corollary}
Let $X: \Sigma^n\rightarrow \mathbb{B}^{n+1}(1)$ be a compact 
connected embedded smooth minimal (or CMC) hypersurface with free boundary  in $S^n(1)$.  
If $\partial \Sigma$ is an isoparametric  hypersurface in $S^n(1)$ which is 
invariant under some compact connected Lie subgroup $G$ of $SO(n+1)$, then $\Sigma$ is $G$-invariant.  
\end{corollary}

Since most of our work takes place locally and an immersion is locally an embedding, we can derive 

\begin{theorem}\label{immersion}
Let $X: \Sigma \to\mathbb{R}^{n+1}$ be an immersion of a connected  hypersurface $\Sigma$ with boundary $\partial\Sigma$ 
such that $X|_{\partial\Sigma}$ is an embedding.
Suppose the same conditions hold as in Theorem \ref{gen2} for 
 a compact connected Lie subgroup $G$ of $SO(n+1)$.  
Then the interior of $X(\Sigma)$ is locally $G$-invariant.

If in addition $X(\Sigma)$ is a closed subset in $\mathbb{R}^{n+1}$ and $X(\partial\Sigma)$ is a $G$-invariant subset in $\mathbb{R}^{n+1}$, then $X(\Sigma)$ is a $G$-invariant set and $X_{\#}[[\Sigma]]$ is a G-invariant current.
\end{theorem}

This paper is organized as follows. Section \ref{sec2} contains some basic definitions and results on group action and real analyticity. In Section \ref{sec3}, we are concerned with CMC hypersurfaces with rotationally symmetric boundary at first and arrive at a preliminary answer to Problem \ref{P1}. Beyond this, we consider isoparametric hypersurfaces in unit spheres to be boundaries in our setting. Based on isoparametric foliations, examples of minimal hypersurfaces and CMC hypersurfaces with $G$-invariant boundaries and $G$-invariant conormals are constructed in Proposition \ref{Prop:construction}.
Section \ref{sec4} is devoted to certain real analytic properties through Morrey's regularity theory. Section \ref{sec5} exhibits proofs of our main results (Theorem \ref{gen2}, case (i)-(iii) and Theorem \ref{localthm}) by making use of the Cauchy-Kovalevskaya theorem. In Section \ref{sec6}, we introduce the definition of Helfrich-type hypersurface and study its real analyticity and symmetry (Theorem \ref{Analyticity2} and Theorem \ref{gen2}, case (iv)). In Section \ref{sec7}, we explore immersion situations and inheritance of symmetry from boundary. For completeness we review the Cauchy-Kovalevskaya theorem and adapt it to our setting in Appendix \ref{secA}.

\section{Preliminaries}\label{sec2}
\subsection{Infinitesimal action, local action and global action}
\begin{definition}[cf. \cite{CG86}]
Let $G$ be a Lie group with Lie algebra $\mathcal{G}$ and $\Sigma$ a smooth manifold. An \emph{infinitesimal action} of $\mathcal{G}$ on $\Sigma$ is a homomorphism of $\mathcal{G}$ to the Lie algebra of smooth vector fields on $\Sigma$. 
A \emph{partial action} $A$ of $G$ on $\Sigma$ is a smooth map 
\begin{equation*}
    \begin{aligned}
    A: \mathcal{D}&\to \Sigma\\
    (g,p)&\to g\cdot p
    \end{aligned}
\end{equation*}
in a neighborhood $\mathcal{D}\subset G \times \Sigma$ of $\{e\}\times\Sigma$ such that $e\cdot p=p$ and $g_1\cdot(g_2\cdot p)=(g_1g_2)\cdot p$ whenever $g_1\cdot(g_2\cdot p)$ lies in $\mathcal{D}$. A partial action defined on $\mathcal{D}=G\times \Sigma$ is called a \emph{global action.}
Two partial actions $(A_1,\mathcal{D}_1)$, $(A_2,\mathcal{D}_2)$ are said to be \emph{equivalent} if there is a domain $\mathcal{D}\subset\mathcal{D}_1\cap \mathcal{D}_2$ containing $\{e\}\times \Sigma$ such that $A_1|_\mathcal{D}=A_2|_\mathcal{D}$. 
A \emph{local action} is an equivalence class of partial actions.
\end{definition}

\begin{remark}
\begin{itemize}
\item[(i).] A global action $A$ defines a local action $A_{loc}$ in an obvious way.
\item[(ii).] The category of local actions is equivalent to the category of infinitesimal actions.
\end{itemize}
\end{remark}

\begin{definition}\label{def2.3}
A subset $\Gamma\subset\Sigma$ is called \emph{locally $A$-invariant} (or \emph{locally $G$-invariant}) if for some representative $(A,\mathcal{D})$ we have $g\cdot p\in\Gamma$ for all $p\in\Gamma$ and $(g,p)\in\mathcal{D}$. A function $f:\Sigma\to\mathbb{R}$ is called \emph{locally $A$-invariant} (or \emph{locally $G$-invariant}) if for any equivalent infinitesimal action $V$ on $\Sigma$ we have $V f\equiv 0$.
\end{definition}

In this paper, we mainly consider Lie subgroups of $SO(n+1)$ and the induced actions. The next proposition guarantees assembling local $G$-invariance to a global one.

\begin{proposition}\label{Killing}
Let $G$ be a connected Lie subgroup of $SO(n+1)$ with Lie algebra $\mathcal{G}$. 
Assume that a connected embedded hypersurface $X: \Sigma \to \mathbb{R}^{n+1}$ with boundary is complete with respect to the induced metric, the interior of $\Sigma$ is locally $G$-invariant and $\partial\Sigma$ is $G$-invariant. Then $\Sigma$ is $G$-invariant.
\end{proposition}
\begin{proof}
For any $\phi\in \mathcal{G}\subset so(n+1)$,  $\phi X$ is a Killing vector field on $\Sigma$ with respect to the induced metric.
Since $\partial\Sigma$ is $G$-invariant, we only need to show that 
$\mathrm{Int}\Sigma$ is $G$-invariant. Actually it suffices to show that 
each maximal integral curve of $\phi X$ in $\mathrm{Int}\Sigma$ is defined for all $t\in\mathbb{R}$. We will prove by contradiction.

Assume that there exists an integral curve $\gamma$
of $\phi X$ in $\mathrm{Int}\Sigma$ 
with $\gamma(0)=p$ and maximal domain $(a,b)$, where $-\infty\le a<0<b<+\infty$. 
 Notice that any integral curve of the Killing vector field $\phi X$ has constant speed. 
Let $\{t_i\}\subset(a,b)$ be a sequence such that $\lim_{i\to\infty}t_i=b$. 
It follows that $\{\gamma(t_i)\}$ is a Cauchy sequence in $\Sigma$. Since $\Sigma$ is complete with respect to the induced metric, $\{\gamma(t_i)\}$ converges to a point $q\in \Sigma$. By the $G$-invariance of $\partial\Sigma$, we know that $q$ does not lie in $\partial\Sigma$. 

 Now choose a neighborhood $U\subset \mathrm{Int}\Sigma$ of $q$ and a positive number $\epsilon$ such that the flow $\psi$ of $\phi X$  is defined on $(-\epsilon,\epsilon)\times U$. Pick $t_k>b-\epsilon$ such that $\gamma(t_k)\in U$ and  define $\tilde{\gamma}:(a,t_k+\epsilon)\to \mathrm{Int}\Sigma$ by
\begin{equation*}
\tilde{\gamma}(t) = \left\{ \begin{array}{ll}
\gamma(t), & a<t<b,\\
\psi_{t-t_k}(\gamma(t_k)), & t_k-\epsilon<t<t_k+\epsilon.
\end{array} \right.
\end{equation*}
$\tilde{\gamma}$ is well-defined since $\psi_{t-t_k}(\gamma(t_k))=\psi_t(p)=\gamma(t)$ for $t<b$. Thus $\tilde{\gamma}$ is an integral curve extending  $\gamma$, which contradicts the maximality of $(a,b)$. The argument for the case $-\infty < a$ is similar. The proof is now complete.
\end{proof}

\subsection{Real analyticity}
We next review  real analyticity of embedded submanifolds in $\mathbb{R}^{n+1}$ (cf. \cite{KP}).
\begin{definition}\label{def4.1}
A subset $\Sigma\subset\mathbb{R}^{n+1}$ is called an $m$-dimensional \emph{real analytic submanifold} if, for each $p\in\Sigma$, there exists an open subset $U\subset\mathbb{R}^{m}$ and a real analytic map $f:U\to\mathbb{R}^{n+1}$ which maps open subsets of $U$ onto relatively open subsets of $\Sigma$ such that
$$
p\in f(U)\quad \textrm{and}\quad \textrm{rank}[Df(u)]=m, \,\forall u\in U,
$$
where $Df(u)$ is the Jacobian matrix of $f$ at $u$.
\end{definition}
\begin{definition}
Let $\mathbb{H}^{m}$ be the closed upper half space of $\mathbb{R}^{m}$. A subset $\Sigma\subset\mathbb{R}^{n+1}$ is called an $m$-dimensional \emph{real analytic submanifold with boundary} if, for each $p\in\Sigma$, there exists an open subset $U\subset\mathbb{H}^{m}$ and a real analytic map $f:U\to\mathbb{R}^{n+1}$ which maps open subsets of $U$ onto relatively open subsets of $\Sigma$ such that
$$
p\in f(U)\quad \textrm{and}\quad \textrm{rank}[Df(u)]=m, \,\forall u\in U,
$$
where $Df(u)$ is the Jacobian matrix of $f$ at $u$.
The pair $(U, f)$ is called a \emph{local parametrization} around $p$.
\end{definition}

Now real analytic functions can be defined on a real analytic submanifold as follows.
\begin{definition}
With the notations above, let $\Sigma$ be a real analytic submanifold (with or without boundary) in $\mathbb{R}^{n+1}$. A function $h:\Sigma\to\mathbb{R}$ is said to be real analytic at $p\in\Sigma$ if there exists a local parametrization $(U, f)$ around $p\in\Sigma$ with $f(q)=p$ such that $h\circ f:U\to\mathbb{R}$ is real analytic at $q\in U$.
\end{definition}

\section{Various examples}\label{sec3}
\subsection{Rotationally symmetric case}
Let $E_1,\cdots, E_{n+1}$ denote the standard basis of $\mathbb{R}^{n+1}$. We start with embedded hypersurface $\Sigma^n\subset\mathbb{R}^{n+1}$ with rotationally symmetric boundary for arbitrary $n\ge2$. Here we allow the boundary $\partial\Sigma$ to be non-connected and only assume that some connected component $\Gamma$ of $\partial\Sigma$ is an $(n-1)$-sphere. 
If the contact angle is constant along $\Gamma$, we obtain 
\begin{proposition}\label{prop3.1}
Let $X:\Sigma^n\to\mathbb{R}^{n+1}$ be an embedding of a connected $C^1$ CMC hypersurface 
 with boundary $\partial\Sigma$ and $\Gamma\subset\partial\Sigma$ be a connected component of $\partial\Sigma$.
Assume that $\Gamma\subset\{x_{n+1}=0\}$ is an $(n-1)$-sphere with constant contact angle $\theta$, then the interior of $\Sigma^n$ is locally rotationally symmetric around $E_{n+1}$.

If in addition $\Sigma$ is complete with respect to the induced metric and $\partial\Sigma$ is rotationally symmetric around $E_{n+1}$, then $\Sigma^n$ is rotationally symmetric around $E_{n+1}$.
\end{proposition}

\begin{remark}
The above result  generalizes  Proposition 5.1 in \cite{PP21} for the case of $n=2$. 
It is a special case of Theorem \ref{gen2} so we omit the proof. 
In this paper we will mainly consider  hypersurfaces with  
general symmetries.
\end{remark}
Note that rotationally symmetric CMC hypersurfaces in $\mathbb{R}^{n+1}$ have been classified by Hsiang and Yu \cite{HY81}. Those hypersurfaces provide natural examples of CMC hypersurfaces with constant contact angle along boundaries.

\subsection{Isoparametric hypersurfaces as boundaries}
A large family of intriguing examples can be constructed based on isoparametric hypersurfaces in spheres. Let us recall some basics about them.

Isoparametric hypersurfaces in $S^{n}(1)$  are compact oriented embedded hypersurfaces with constant principal curvatures. 
By a fundamental result  due to M{\"u}nzner \cite{Mun80},  the number $g$ of distinct principal curvatures $\lambda_1>\lambda_2 >\cdots >\lambda_g$
of an isoparametric hypersurface in $S^n(1)$ must be $1, 2, 3, 4$ or $6$ and their multiplicities satisfy $m_i=m_{i+2}$ with index modulo $g$.  
Any isoparametric hypersurface in $S^n(1)$  is given by a regular level sets of the restriction of a homogeneous polynomial $F: \mathbb{R}^{n+1} \rightarrow \mathbb{R}$ of degree $g$ in $\mathbb{R}^{n+1}$ to $S^n(1)$ satisfying the so-called Cartan-M\"{u}zner equations:  
\begin{equation*}
\left\{ \begin{array}{ll}
|\nabla F|^2= g^2|x|^{2g-2}, \\
~~~\triangle F~~=\dfrac{m_2-m_1}{2}g^2|x|^{g-2}.
\end{array}\right.
\end{equation*}
Such a polynomial $F$ is called the \emph{Cartan-M{\"u}nzner polynomial}, and 
$f=F~|_{S^{n}(1)}$ takes values in $[-1, 1]$. 
The level sets $M_+=f^{-1}(1)$ and  $M_-=f^{-1}(-1)$, called the focal submanifolds, are smooth minimal submanifolds of $S^{n}(1)$ with codimensions $m_1+1$ and $m_2+1$ respectively (cf. \cite{Nom73}). 

Denote the principal curvature by $\lambda_i=\cot(\varphi+\frac{i-1}{g}\pi)$ with $\varphi\in (0, \frac{\pi}{g})$ $(i=1,\ldots,g)$.
We collect the following useful properties of isoparametric hypersurfaces:
\begin{itemize}
\item[(i).] $\mathrm{dist}(M_+,M_-)=\frac{\pi}{g}$, where $g\in \{1, 2, 3, 4, 6\}$.

\item[(ii).] For each $\varphi\in (0, \frac{\pi}{g})$, let $M_{\varphi}$ be the tube of constant radius $\varphi$ around $M_+$ in $S^{n}(1)$. Then there exists some $c\in (-1, 1)$, such that $M_{\varphi}=f^{-1}(c)$. 

\item[(iii).] For each $\varphi\in (0, \frac{\pi}{g})$, the volume $V(M_{\varphi})$ of $M_{\varphi}$ is given by $\left(\sin\frac{g\varphi}{2}\right)^{m_1}\left(\cos \frac{g\varphi}{2}\right)^{m_2}$ up to a positive constant independent of $\varphi$, and the mean curvature of $M_{\varphi}$ is $\frac{d}{d\varphi}\log V(M_{\varphi})$ with respect to the unit normal $\frac{\nabla f}{|\nabla f|}$.

\item[(iv).] Each isoparametric hypersurface in $S^{n}(1)$ is actually a real analytic submanifold in $\mathbb{R}^{n+1}$ of codimension $2$, which is an intersection of some level set of Cartan-M{\"u}nzner polynomial $F$ and $S^{n}(1)$.
\end{itemize}

Due to well-developed classification results, isoparametric hypersurfaces in unit spheres are either homogeneous or of OT-FKM type with $g=4$. 
The latter can provide infinitely many inhomogeneous isoparametric hypersurfaces. For the classification of isoparametric hypersurfaces in spheres and related topics,  we refer to \cite{Qi, QT14}.

\subsubsection{Homogeneous isoparametric hypersurfaces}
A hypersurface $M^{n-1}$ in $S^{n}(1)$ is homogeneous if 
$M^{n-1}$ is an orbit of a closed connected Lie subgroup of $G\subset SO(n+1)$.
By virtue of Hsiang-Lawson \cite{HL71} and Takagi-Takahashi \cite{TT72}, 
any homogeneous isoparametric hypersurface $M^{n-1}$ in $S^{n}(1)$ can be obtained as a principal orbit of a linear isotropy representation of compact Riemannian symmetric pairs $(U,G)$ of rank $2$. 

As a corollary of Theorem \ref{gen2} and the property (iv), when the boundary is a homogeneous isoparametric hypersurface, we obtain
\begin{corollary}\label{corIso}
Let $\Sigma^n$ be an embedded hypersurface in $\mathbb{R}^{n+1}$ with boundary $\partial\Sigma$. Assume that one of the connected components of $\partial\Sigma$ is a homogeneous isoparametric hypersurface $M^{n-1}$ in $S^{n}(1)$ as an $G$-orbit. 
If the contact angle is constant along $M^{n-1}$ and $\Sigma$ satisfies the same condition as in Theorem \ref{gen2},
then the interior of $\Sigma$ is locally $G$-invariant.

If in addition $\Sigma$ is complete with respect to the induced metric, and each connected component of $\partial\Sigma$ is a $G$-invariant submanifold in $\mathbb{R}^{n+1}$, then $\Sigma$ is $G$-invariant.
\end{corollary} 

\subsubsection{Isoparametric hypersurfaces of OT-FKM type}
For a given symmetric Clifford system $\{P_0, P_1,\cdots, P_{m}\}$ on $\mathbb{R}^{2l}$ satisfying
$P_{\alpha}P_{\beta}+P_{\beta}P_{\alpha}=2\delta_{\alpha\beta}Id$\,\,
for $0 \leq \alpha,\beta\leq m$,
Ferus, Karcher and
M\"{u}nzner \cite{FKM81} constructed a Cartan-M\"{u}nzner polynomial $F$ of degree $4$ on
$\mathbb{R}^{2l}$
\begin{eqnarray*}\label{FKM isop. poly.}
&&\qquad F:\quad \mathbb{R}^{2l}\rightarrow \mathbb{R}\nonumber\\
&&F(x) = |x|^4 - 2\displaystyle\sum_{\alpha = 0}^{m}{\langle
P_{\alpha}x,x\rangle^2}.
\end{eqnarray*}
A regular level set of $f=F|_{S^{2l-1}(1)}$ is called an isoparametric hypersurface of \emph{OT-FKM type}. Its multiplicity pair is $(m_1, m_2)=(m, l-m-1)$  provided $m>0$ and $l-m-1>0$, where $l=k\delta(m)$, $k$ being a positive integer and $\delta(m)$ the dimension of irreducible module of the Clifford algebra $\mathcal{C}_{m-1}$.

In general, isoparametric hypersurfaces of OT-FKM type are not extrinsic homogeneous. 
Instead, we have
\begin{proposition}[\cite{FKM81}]\label{spin action}
Let $M^{2l-2}$ be an isoparametric hypersurface of OT-FKM type in $S^{2l-1}(1)$ with $(m_1, m_2)=(m, l-m-1)$, and $Spin(m+1)$ be the connected Lie subgroup in $SO(2l)$ generated by the Lie subalgebra $\mathrm{Span }\{P_{\alpha}P_{\beta}~|~0\leq \alpha, \beta\leq m\}\subset \mathfrak{so}(2l)$. Then $M^{2l-2}$ is $Spin(m+1)$-invariant.
\end{proposition}
\begin{remark} As a simple illustration of the $Spin(m+1)$ action in Proposition \ref{spin action}, a concrete example is given as follows.
For $m=1$, $l\ge3$, the corresponding
isoparametric hypersurface is diffeomorphic to $S^1\times V_{2}(\mathbb{R}^l)$, where $V_{2}(\mathbb{R}^l)=\{(x,y)~|~x,y\in \mathbb{R}^l, |x|=|y|=1, x\bot y\}$. Now the $Spin(2)$ action is given by 
\begin{equation*}
\begin{aligned}
Spin(2)&\times (S^1\times V_{2}(\mathbb{R}^l))&\longrightarrow&\quad\quad\quad S^1\times V_{2}(\mathbb{R}^l),\\
(e^{it}&, (e^{i s}, (x,y)))&\longmapsto&\quad (e^{i( s-2t)},(\cos{t}x+\sin{t}y, -\sin{t}x+\cos{t}y)).
\end{aligned}
\end{equation*}
\end{remark}

When the boundary is an isoparametric hypersurface of OT-FKM type, we have the following  application of the property (iv), Proposition \ref{spin action} and Theorem \ref{gen2}.

\begin{corollary}\label{OT-FKM}
Let $\Sigma^{2l-1}$ be an embedded hypersurface in $\mathbb{R}^{2l}$ with boundary $\partial\Sigma$. Assume that a connected component of $\partial\Sigma$ is an isoparametric hypersurface $M^{2l-2}$ of OT-FKM type in $S^{2l-1}(1)$ which is $Spin(m+1)$-invariant. Moreover, assume that $\Sigma$ satisfies the same condition as in Theorem \ref{gen2} with $G=Spin(m+1)$. Then the interior of $\Sigma$ is locally $Spin(m+1)$-invariant.

If in addition $\Sigma$ is complete with respect to the induced metric, and $\partial\Sigma$ is $Spin(m+1)$-invariant, then $\Sigma$ is $Spin(m+1)$-invariant.
\end{corollary}
         
\subsection{Constructions via isoparametric foliations}  
Now we employ isoparametric foliations to  construct examples of foliated minimal hypersurfaces and CMC hypersurfaces with isoparametric boundaries and certain prescribed conormals.

An immersed hypersurface $X:\Sigma^n\to\mathbb{R}^{n+1}$ is called $G$-invariant if there exists a smooth $G$-action on $\Sigma$ such that $X\circ h=h\circ X$ for any $h\in G$.         

Recall that a $G$-invariant immersed hypersurface $\Sigma$ is minimal if and only if $\Sigma$ is minimal among $G$-invariant competitors (see Hsiang-Lawson \cite{HL71}).
Therefore, when the $G$-action on $\mathbb{R}^{n+1}$ is of cohomogeneity $2$, 
there is a one to one correspondence between $G$-invariant minimal hypersurfaces $\Sigma$ and  geodesics in the orbit space $\mathbb{R}^{n+1}/G$.
Given a $G$-homogeneous isoparametric boundary $M$ and any prescribed $G$-invariant outward unit conormal $V$ along $M$, we investigate the geodesic in the orbit space 
                    $$
                    \mathscr R= C\left(\frac{\pi}{g}\right):=
                      \left\{ re^{i\theta} \in \R^2: 0\leq \theta \leq \frac{\pi}{g}  \right\}
                    $$
               with the canonical metric
                  $$
                  g_c=\underline V^2 \left(dx^2+dy^2\right),
                  $$
                  starting from $p_M=M\slash G$ in $\mathscr R$. 
                  Here $g$ is the number of distinct principal curvatures of 
                  $M$ and $\underline V=\underline V(x,y)$ is the volume of $G$-invariant orbit in $\R^{n+1}$ represented by the  point $(x,y)$.

\begin{figure}[htbp]
	\centering
	\begin{subfigure}[t]{0.43\textwidth}
		\centering
		\includegraphics[scale=0.55]{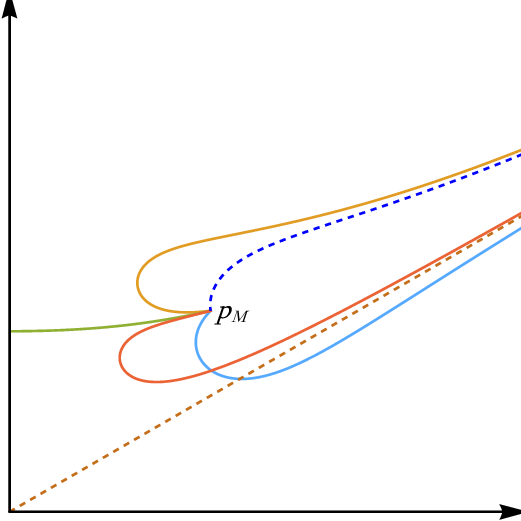}
                              \captionsetup{font={scriptsize}} 
                               \caption{$g=2$, $m_1=2$, $m_2=6$ and $p_M$ is not on the minimal cone $\varphi=\frac{\pi}{6}$}
                               \label{fig:1a}
	\end{subfigure}
	\quad\quad \quad\quad 
	\begin{subfigure}[t]{0.43\textwidth}
		\centering
	 \includegraphics[scale=0.55]{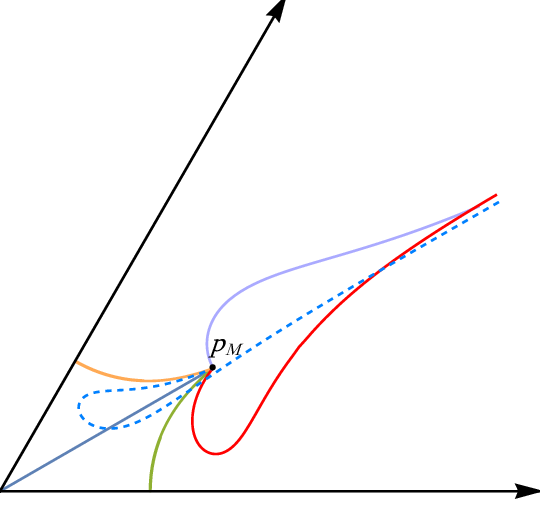}
                           \captionsetup{font={scriptsize}} 
                             \caption{$g=3$, $m_1=m_2=2$}
                             \label{fig:1b}
	\end{subfigure}
	\\
	\begin{subfigure}[t]{0.43\textwidth}
		\centering
		\includegraphics[scale=0.55]{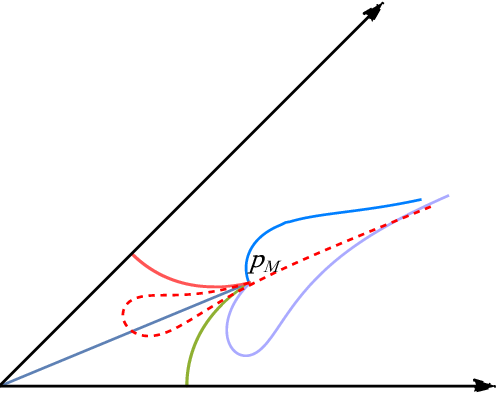}
                              \captionsetup{font={scriptsize}} 
                               \caption{$g=4$, $m_1=m_2=2$}
                               \label{fig:1c}
	\end{subfigure}
	\quad\quad \quad\quad 
	\begin{subfigure}[t]{0.43\textwidth}
		\centering
	 \includegraphics[scale=0.55]{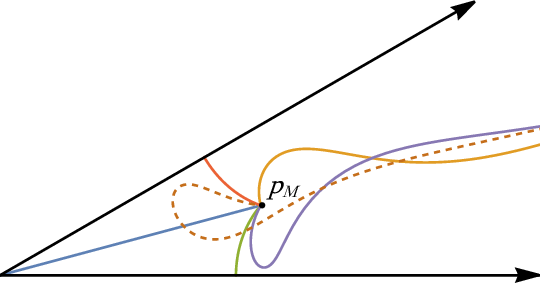}
                           \captionsetup{font={scriptsize}} 
                             \caption{$g=6$, $m_1=m_2=1$}
                             \label{fig:1d}
	\end{subfigure}
	\captionsetup{font={small}} 
	\caption{Complete geodesics starting from $p_M$ in orbit space $\mathscr{R}$}\label{fig:1}
\end{figure}

              The geometric behavior of the complete geodesics in $\mathscr R$  was studied by Wang in \cite{Wang94}. Note that a complete geodesic $\gamma: J\longrightarrow\mathscr R$  starting from $p_M$ 
              where $J=[0,1]$ or $[0,\infty)$ is determined by the direction $-\gamma'(0)$ that corresponds to $V$.  According to Lemma 3.1 and Theorem 4.8 in \cite{Wang94}, there are three different types of $\gamma$ depending on the choice of $p_M$ and $V$: 
(i) $\gamma$ goes to the origin $O$ forming a minimal cone; (ii) $\gamma$ hits $\partial \mathscr R$ perpendicularly; (iii) $\gamma$ extends to infinity asymptotically toward the minimal cone (see Fig. \ref{fig:1}). Moreover, among all geodesics starting from $p_M$, there are at most one of them goes to the origin $O$ and at most two of them hit $\partial \mathscr R\backslash \{O\}$. We would like to remark that there is no geodesic with both ending points in $\p \mathscr R$
              since there is no closed minimal hypersurface in $\R^{n+1}$.

              In fact, one can obtain the same type results for non-homogeneous isoparametric hypersurface $M$ of $S^n(1)$.
               Note that $\R^{n+1}\backslash \{O\}$ is foliated by positive homotheties of leaves of an isoparametric foliation of $S^n(1)$.
               Now we define an ``orbit space" $\mathscr R=C\left(\frac{\pi}{g}\right)$
               \ in the same way endowed with metric $g_c$ where, up to a positive constant,
               \begin{equation}\label{inh}
               g_c=
               {\underline V}^2
               \left(
               dr^2+ r^2 d\varphi^2
               \right)
               \end{equation}
               in the polar coordinate, where
               $$
               \underline V = 
               r^{n-1} \cdot 
               \left(\sin \frac{g\varphi}{2}\right)^{m_1} 
               \left(\cos \frac{g\varphi}{2}\right)^{m_2}.
               $$

According to Ferus and Karcher \cite{FK85}, for any immersed curve $\gamma(s)=(r(s),\varphi(s))$ in $\mathscr R$ parametrized by its arc length, the mean curvature $\tilde h$ of the hypersurface $\Sigma$ (in $\R^{n+1}$)
                      corresponding to $\gamma$ is determined by
               \begin{equation}\label{isopH}
               \left\{\begin{array}{l}
               r'(s)=\sin\alpha,\\
               \varphi'(s)=\frac{\cos\alpha}{r},\\
               \alpha'(s)=-\tilde{h}+n\frac{\cos\alpha}{r}-h(\varphi)\frac{\sin\alpha}{r},
               \end{array}\right.
               \end{equation}
               where $h(\varphi)=\frac{g}{2}(m_1\cot{\frac{g\varphi}{2}}-m_2\tan{\frac{g\varphi}{2}})$ is the mean curvature (with respect to the unit normal toward $M_+$)
               of the isoparametric hypersurface of distance $\varphi$ to $M_+$ in $S^n(1)$ and $\alpha$ is the angle between $\gamma'(s)$ and $\frac{\partial}{\partial\varphi}$. By Definition \ref{defAngle}, the contact angle of $\Sigma$ along $M$ is $\frac{\pi}{2}-\alpha(0) \textrm{ mod } 2\pi$.
               
              \begin{figure}[htbp]
	\centering
	\begin{subfigure}[t]{0.43\textwidth}
		\centering
		\includegraphics[scale=0.6]{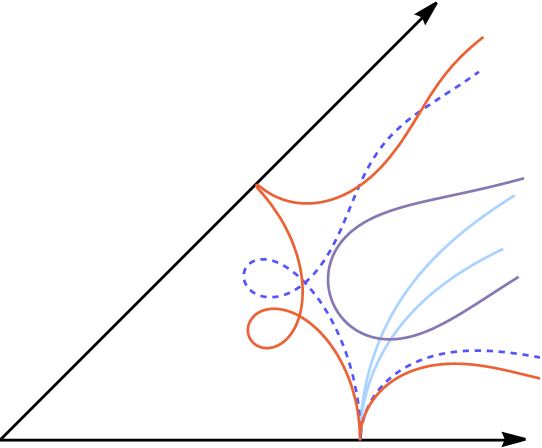}
                              \captionsetup{font={scriptsize}} 
                               \caption{Global solution curves with different $\tilde{h}$}
                               \label{fig:2a}
	\end{subfigure}
	\quad\quad \quad\quad 
	\begin{subfigure}[t]{0.43\textwidth}
		\centering
	 \includegraphics[scale=0.55]{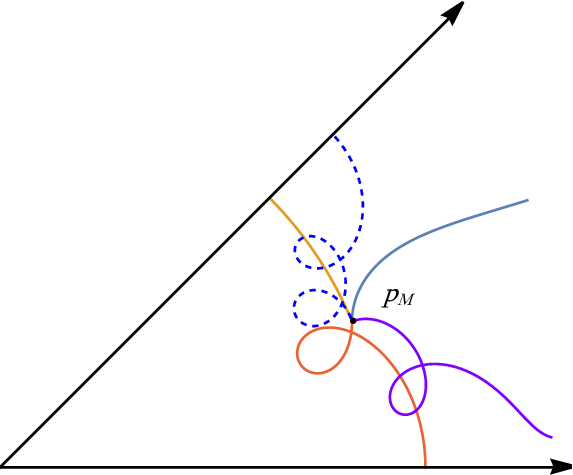}
                           \captionsetup{font={scriptsize}} 
                             \caption{Solution curves starting from $p_M$ with different $\alpha$ and $\tilde{h}$}
                             \label{fig:2b}
	\end{subfigure}
	\\
\begin{subfigure}[t]{0.43\textwidth}
\centering
\includegraphics[scale=0.55]{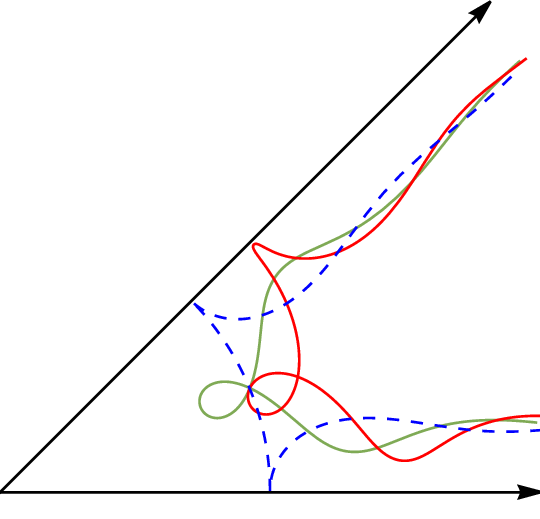}
\captionsetup{font={scriptsize}} 
\caption{Global solution curves with same $\tilde{h}$}
\label{fig:2c}
\end{subfigure}
	\quad\quad \quad\quad 
	\begin{subfigure}[t]{0.43\textwidth}
		\centering
	 \includegraphics[scale=0.5]{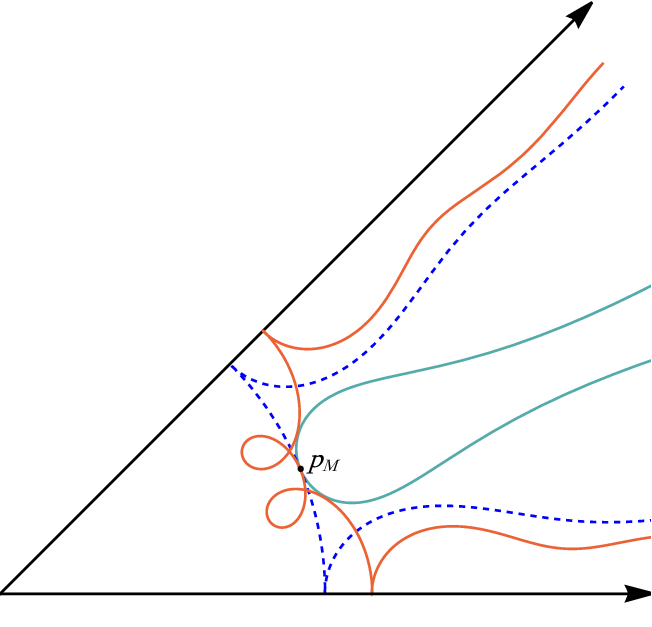}
                           \captionsetup{font={scriptsize}} 
                             \caption{Solution curves starting from $p_M$ with $\alpha=0$ and different $\tilde{h}$}
                             \label{fig:2d}
	\end{subfigure}
	\captionsetup{font={small}} 
	\caption{Solution curves of \eqref{isopH} in ``orbit space" $\mathscr R$  with $g=4$, $m_1=m_2=2$}\label{fig:2}
\end{figure}
             
              Therefore, given homogeneous or non-homogeneous isoparametric boundary $M$
              and prescribed constant contact angle along $M$, the existence question for hypersurface of non-zero constant mean curvature $\tilde h$ can be solved by the ODE system \eqref{isopH}. 
              We notice that the geometric behavior of solution curves of \eqref{isopH} was studied by Hsiang and Huynh in \cite{HH87}. On the one hand, Theorem A in \cite {HH87} states that every global solution curve has two asymptotic lines parallel to the boundary lines $\varphi=0$ and $\varphi=\pi/g$ respectively. On the other hand, Theorem B in \cite{HH87} says that a global solution curve can hit the boundary line $\varphi=0$ (resp. $\varphi=\pi/g$) at most once. As a result, a solution curve $\gamma$ of \eqref{isopH} starting from $p_M$ either hits $\partial\mathscr{R}$ (in fact, perpendicularly) or extends to infinity. Moreover, either case produces both properly immersed and embedded examples (see Fig. \ref{fig:2}). 
              
              Apparently, an $\tilde h$-solution curve generates a corresponding smooth $\tilde h$-CMC hypersurface.
              In particular, a $0$-solution curve for
              minimal hypersurface 
              is indeed a geodesic with respect to \eqref{inh}. In fact, ${\tilde{h}}/{\underline{V}}$ is exactly the geodesic curvature of $\gamma$ with respect to $g_c$ \eqref{inh}.
              
Based on the discussion above, we summarize this section by the following
\begin{proposition}\label{Prop:construction}
Given any isoparametric hypersurface $M^{n-1}$ in $S^{n}(1)$ and any angle $\varphi\in[0, 2\pi)$, it holds
\begin{itemize}
\item[\rm{(i)}] \textbf{The minimal case:} Either there exists a minimal cone $\Sigma$ in $\mathbb{R}^{n+1}$ or a complete properly immersed 
minimal hypersurface $X: \Sigma\to \mathbb{R}^{n+1}$ with $X(\partial \Sigma)=M$ and constant contact angle $\varphi$ along $M$.
\item[\rm{(ii)}] \textbf{The CMC case:} For any $\tilde{h}\ne 0$, there exists a complete properly immersed  hypersurface $X:\Sigma\to \mathbb{R}^{n+1}$ of constant mean curvature $\tilde{h}$ with $X(\partial \Sigma)=M$ and constant contact angle $\varphi$ along $M$.
\end{itemize}
\end{proposition}

\begin{remark} 
 With the fixed boundary $M^{n-1}$, in case (i) there is at least one compact embedded minimal hypersurface;
 in case (ii) 
 there are infinitely many geometrically distinct compact immersed $\tilde{h}$-CMC hypersurfaces in $\mathbb{R}^{n+1}$. 
\end{remark}

\begin{remark}
When $\Sigma$ is connected we have the uniqueness in the language of integral current, i.e., in either case (i) or (ii) every immersion map $X$ with $X(\partial \Sigma)=M$ and constant contact angle $\varphi$ along $M$ in the statement induces the same integral current $X_\#[[\Sigma]]$ (cf. \cite{FF}).
\end{remark}

\section{Real analyticity of hypersurfaces in \texorpdfstring{$\mathbb{R}^{n+1}$}{Euclidean space} }\label{sec4}

Next we review materials on real analyticity, and prove Theorem \ref{Analyticity} and Theorem \ref{Analyticity2}, which will play an important role in this paper.

\subsection{Morrey's regularity theory}
To loosen the restriction on the regularity requirement of  hypersurfaces considered in our main results,
we recall the interior and boundary regularity theory of Morrey \cite{Morrey58-1, Morrey58-2}.
For the reader's convenience, we will introduce some basic notations before stating Morrey's results.

Let $\mathcal{D}\subset \mathbb{R}^p$ be a bounded domain, $\bar{\mathcal{D}}$ be its closure and $x=(x_1,\cdots,x_p)\in \mathcal{D}$. Consider a system of nonlinear partial differential equations in $u=(u^1,\cdots, u^N): \bar{\mathcal{D}}\to \mathbb{R}^N$,
\begin{equation}\label{nLS}
\phi_j(x,u,Du,D^2u,\cdots)=0,\quad j=1,\cdots,N.
\end{equation}
This system is called a \emph{real analytic system} if each $\phi_j$ is real analytic for all values of its arguments. The linearization of the nonlinear system \eqref{nLS} along $u$ is defined by
\begin{equation}\label{LS}
\begin{aligned}
 L_{jk}(x, D)v^k:=\frac{d}{d \rho}\phi_j (x, u+\rho v,D(u+\rho v),D^{2}(u+\rho v),\cdots){\Big|}_{\rho=0}=0,&\\ v=(v^1,\cdots,v^N),\quad j,k=1,\cdots, N.&
\end{aligned}
\end{equation}

\subsubsection{Interior regularity}\label{4.2.1}
Assume that there exist integers $s_1,\cdots,s_N$ and $t_1,\cdots,t_N$ such that each operator $L_{jk}(x,D)$ is of order not greater than $s_j+t_k$. Let $L^0_{jk}(x,D)$ be the terms in $L_{jk}(x,D)$ which are exactly of order $s_j+t_k$ and denote its characteristic polynomial by $L^0_{jk}(x,\lambda)$, where $\lambda=(\lambda_1,\cdots,\lambda_p)$.
\begin{definition}[cf. \cite{ND55}]
The linear system \eqref{LS} is called \emph{elliptic} if such integers $s_1,\cdots,s_N$ and $t_1,\cdots,t_N$ exist that at each point $x$ the determinant $L(x,\lambda):=det(L^0_{jk}(x,\lambda))$ of the characteristic polynomial is not zero for any non-zero $\lambda\in\mathbb{R}^p$.
\end{definition}
\begin{definition}[cf. \cite{ND55}]
The nonlinear system \eqref{nLS} is called \emph{elliptic along $u$} if its linearization \eqref{LS} along $u$ form a linear elliptic system.
\end{definition}
Since we can add the same integer to all the $t_j$ and subtract it from all the $s_j$, we assume that $max\{s_j |1\le j \le N\}=0$ in Section \ref{4.2.1}. Now we are ready to state the following theorem.
\begin{theorem}[cf. \cite{Morrey58-1}]\label{Morrey1}
Let $\mathcal{D}\subset \mathbb{R}^p$ be a bounded domain and $\bar{\mathcal{D}}$ be its closure. Let $u=(u^1,\cdots,u^N): \bar{\mathcal{D}}\to \mathbb{R}^N$ be a solution of a real analytic nonlinear system
\begin{equation}\label{AES}
\phi_j(x_1,\cdots,x_p,\cdots,u^k,D u^k,\cdots,D^{s_j+t_k}u^k,\cdots)=0,\quad j=1,\cdots,N,
\end{equation}
where $\phi_j$ involves derivatives of $u^k$ of order not greater than $s_j+t_k$ for $1\le j,k\le N$ and $max\{s_j |1\le j \le N\}=0$.
Assume that the system \eqref{AES} is elliptic along $u$ and each $u^k$ is of class $C^{t_k, \mu}$ in $\bar{\mathcal{D}}$ with some $0<\mu<1$, then $u$ is real analytic at each interior point of $\mathcal{D}$.
\end{theorem}

\subsubsection{Boundary regularity}
Assume that there exist non-negative integers $s_1,\cdots,s_N$ such that each operator $L_{jk}(x,D)$ is of order not greater than $s_j+s_k$ and define $s:=max\{s_j|1\le j\le N\}$. Let $L^0_{jk}(x,D)$ be the terms in $L_{jk}(x,D)$ which are exactly of order $s_j+s_k$ and denote its characteristic polynomial by $L^0_{jk}(x,\lambda)$, where $\lambda=(\lambda_1,\cdots,\lambda_p)$.
\begin{definition}[cf. \cite{Nirenberg55}]
The linear system \eqref{LS} is called \emph{strongly elliptic} if there exist  integers $s_1,\cdots,s_N$ such that at each point $x$ the characteristic matrix $(L^0_{jk}(x,\lambda))$ is definite, i.e., $L^0_{jk}(x, \lambda)\xi^j\bar{\xi}^k\ne0$ for any non-zero $\lambda\in\mathbb{R}^{p}$ and any non-zero $\xi\in\mathbb{C}^N$.
Moreover, the linear system \eqref{LS} is called \emph{uniformly strongly elliptic along $u$ in $\bar{\mathcal{D}}$} if it is strongly elliptic and 
\begin{equation}
\textrm{Re}[L^0_{jk}(x, \lambda)\xi^j\bar{\xi}^k]\ge M \sum_{j=1}^N |\lambda|^{2s_j}|\xi_j|^2, \quad \forall \lambda\in \mathbb{R}^p, \, \forall \xi\in \mathbb{C}^N,
\end{equation}
for some $M>0$ which is independent of $\lambda$, $\xi$ and $x\in\bar{\mathcal{D}}$.
\end{definition}
\begin{definition}[cf. \cite{Nirenberg55}]
The nonlinear system \eqref{nLS} is called \emph{strongly elliptic (uniformly strongly elliptic, respectively) along $u$} if its linearization \eqref{LS} is strongly elliptic (uniformly strongly elliptic, respectively).
\end{definition}
Then Morrey's boundary regularity theorem can be stated as follows.
\begin{theorem}[cf. \cite{Morrey58-2}]\label{Morrey2}
Let $\mathcal{D}\subset \mathbb{R}^p$ be a bounded domain and $\bar{\mathcal{D}}$ be its closure. Let $u=(u^1,\cdots,u^N): \bar{\mathcal{D}}\to \mathbb{R}^N$ be a solution of the following real analytic nonlinear system 
\begin{equation}\label{ASES}
\phi_j(x_1,\cdots,x_p,\cdots,u^k,D u^k,\cdots,D^{s_j+s_k}u^k,\cdots)=0,\quad j=1,\cdots,N,
\end{equation}
where $\phi_j$ involves derivatives of $u^k$ of order not greater than $s_j+s_k$ for $1\le j,k\le N$. Assume that the system \eqref{ASES} is uniformly strongly elliptic along $u$ in $\bar{\mathcal{D}}$ and each $u^k$ is of class $C^{s+s_k,\mu}$ in $\bar{\mathcal{D}}$ with some $0<\mu<1$ and $s:=max\{s_j|1\le j\le N\}$. If $u$ possesses Dirichlet data which are real analytic along a real analytic portion $\mathcal{C}$ of the boundary of $\bar{\mathcal{D}}$, then $u$ can be extended real analytically across  $\mathcal{C}$.
\end{theorem}

\subsection{Real analyticity of hypersurfaces}
As an application of Morrey's regularity theorems, we get the following
\begin{theorem}\label{Analyticity}
Let $X:\Sigma \rightarrow \mathbb{R}^{n+1}$  be  an embedded hypersurface with boundary $\partial\Sigma$. If $\Sigma$ is one of the first three cases listed in Theorem \ref{gen2},
 then $\Sigma$ is real analytic at each interior point of $\Sigma$ and can be extended analytically across each real analytic portion of the boundary $\partial\Sigma$.
\end{theorem}
\begin{proof} We divide the proof into three parts. In the first two parts we deal with cases (i) and (ii) simultaneously, and in Part 3 we deal with case (iii).

\textbf{Part 1:} Firstly, we consider minimal hypersurfaces and hypersurfaces of non-zero constant mean curvature as critical hypersurfaces of the functionals $A(\Sigma)$ and $J(\Sigma)$, respectively. Here $A(\Sigma):=Area(\Sigma)$ and $J(\Sigma):=A(\Sigma)+\frac{n}{n+1} H_0\int_{\Sigma}\langle X, \nu\rangle d\sigma$, where $H_0:=A(\Sigma)^{-1}\int_{\Sigma}H\textrm{d}A$, $\nu$ is a unit normal vector field and $H$ is the mean curvature with respect to $\nu$. According to Theorem 9.2 in \cite{Morrey54},  it follows that  $C^1$ minimal hypersurfaces and $C^1$ CMC hypersurfaces are of class $C^{2,\alpha}$.

\textbf{Part 2:} Secondly, since minimal hypersurfaces have vanishing  
mean curvature, we only need to deal with $C^{2,\alpha}$ CMC hypersurfaces.

For each interior point $p$ in $\Sigma$, without loss of generality we 
choose an orthonormal basis $\{e_A\}_{1\le A\le n+1}$ of $\mathbb{R}^{n+1}$ such that $p$ is the origin and $\nu(p)=e_{n+1}$. Let $\pi: \mathbb{R}^{n+1}\to T_p\Sigma$ be the orthogonal projection to $T_p\Sigma$. Then there is a neighborhood $V$ of $p$ such that $V\cap\Sigma$ can be regarded as a graph defined on $\pi(V\cap\Sigma)\subset T_p\Sigma$, i.e., there exists $u:\pi(V\cap\Sigma)\to\mathbb{R}$ such that $X(q)=(x, u(x))$ for any $q\in V\cap\Sigma$ and $x=\pi(q)$. By Part 1, $u$ is locally $C^{2,\alpha}$ around $p$.

Under this graph representation, it is well known that $u$ satisfies the following quasilinear elliptic equation of second order:
\begin{equation}\label{CMCeq}
\sum_{i,j=1}^n(\delta_{ij}-\frac{u_i u_j}{W^2})u_{ij}=  n HW,
\end{equation}
where $u_i=\frac{\partial}{\partial x_i}u$ and $W=(1+\sum_i u_i^2)^{\frac{1}{2}}$.
Thus as a $C^{2,\alpha}$ solution to the quasilinear real analytic elliptic equation \eqref{CMCeq}, $u$ is real analytic at each interior point by Theorem \ref{Morrey1}.

As for the boundary part, for each point $p$ on a real analytic portion $\Gamma\subset\partial\Sigma$, we also have the above graph representation. Moreover, we can choose a local parametrization $(U,f)$ of $p$ such that $f(U)=V\cap\Gamma$, where $f:U\subset\mathbb{R}^{n-1}\to\mathbb{R}^{n+1}$ is a real analytic map with rank $n-1$.

We need to show that $u$ possesses real analytic Dirichlet data along $\pi(f(U))$. Since $\pi\circ f:U\to T_p\Sigma$ is real analytic with rank $n-1$, $\pi(f(U))$ is a real analytic hypersurface of $T_p\Sigma\cong \mathbb{R}^n$. More precisely, $(\pi(f(U)),\pi\circ f)$ is a real analytic parametrization of $\pi(f(U))$ in $T_p\Sigma$. For each $q\in U$ we have 
$$
u\circ \pi\circ f(q)= \langle f(q), e_{n+1}\rangle,
$$
which is real analytic in $U$. Hence $u$ is real analytic in $\pi(f(U))$. Note that each $u_i$ $(1\le i \le n)$ is uniformly bounded around $p$, we can see that \eqref{CMCeq} is uniformly strongly elliptic along $u$ in some closed neighborhood of $p$.

Therefore, as a $C^{2,\alpha}$ solution to the quasilinear real analytic strongly elliptic equation \eqref{CMCeq}, $u$ can be extended analytically to a neighborhood of $p\in \mathbb{R}^{n+1}$ by Theorem \ref{Morrey2}. This completes the proof of the first two cases.

\textbf{Part 3:} Finally, to deal with the case of constant $r$-th mean curvature, let us first introduce the $L_r$ operator.

Choose a local orthonormal frame $\{e_1,\cdots, e_n\}$ of $\Sigma$ and let $\nu$ be a unit normal vector field. Denote the components of the second fundamental form with respect to $\nu$ by $(h_{ij})$, then the classical Newton transformations $T^r$ are defined inductively by
\begin{equation*}
\begin{aligned}
T^0_{ij}&=\delta_{ij},\\
T^1_{ij}&=S_1\delta_{ij}-h_{ij},\\
&\cdots\\
T^r_{ij}&=S_r\delta_{ij}-\sum_k T^{r-1}_{ik}h_{kj},
\end{aligned}
\end{equation*}
for $i, j=1,\cdots, n$, 
where $S_r$  is the $r$-th elementary symmetric polynomial of the principal curvatures.
Associated to each Newton transformation $T^r$, there is a second-order differential operator $L_r$ defined by
\begin{equation}\label{LrDef}
L_rf:=\sum_{i,j}T^r_{ij}f_{ij},\quad\textrm{for any } f\in C^\infty(\Sigma),
\end{equation}
where $f_{ij}$ is the Hessian of $f$.
When $r=0$, $L_0$ is just the Laplace-Beltrami operator on $\Sigma$. It is also known (cf. \cite{BC97}) that $L_{r-1}$ satisfies
\begin{equation}\label{eq3.4}
L_{r-1} \langle X,a\rangle=rS_r\langle\nu,a\rangle,
\end{equation}
for any fixed vector $a\in\mathbb{R}^{n+1}$.
In light of the proof for Proposition 3.2 in \cite{BC97}, when $H_r>0$ and in addition there exists an interior elliptic point, it follows that $L_{r-1}$ is elliptic for the hypersurface $\Sigma$ with boundary.

Similar to cases (i) and (ii), for each point $p$ in $\Sigma$ we 
regard a neighborhood of $p$ in $\Sigma$ as a graph $X=(x,u(x))$ defined on a domain in $T_p\Sigma$. From the assumption that $u$ is $C^{2}$,  by \eqref{eq3.4} we  
derive that $u$ satisfies the following elliptic equation of second order:
\begin{equation}\label{CrMCeq}
L_{r-1} u= \frac{rS_r}{W} .
\end{equation}

Since \eqref{CrMCeq} is a nonlinear real analytic elliptic equation, applying Theorem 9.1 in \cite{Morrey54} we 
obtain that $u$ is of class $C^{2,\alpha}$. Therefore the conclusion follows from Theorem \ref{Morrey1} and Theorem \ref{Morrey2} as before.
\end{proof}
\begin{remark}
We can also use the regularity theory of second-order elliptic equations to prove \textbf{Part 1} (cf. \cite{GT} Section 7 and Section 9).
Furthermore,  the $C^1$ condition in case (i)  can be weakened to be Lipschitz continuous.
\end{remark}

\section{Proof of the main theorem}\label{sec5}

Now let us consider arbitrary compact connected Lie subgroup $G\subset SO(n+1)$ beyond $SO(n)$.
As mentioned in Remark \ref{rk3},
our method only relies on  local assumptions along boundary and  we can establish the following 

\begin{theorem}\label{localthm}
Let  $X:\Sigma \to\mathbb{R}^{n+1}$ be an embedded hypersurface with boundary $\partial\Sigma$.
Suppose that $U\subset \partial\Sigma\bigcap S^n(R)$ for some positive $R$  is  a real analytic open piece of $\partial\Sigma$,
and moreover, that $U$ is locally $G$-invariant.
If the contact angle $\theta$ along $U$ is locally $G$-invariant and $\Sigma$ is among $\rm{(i)}$-$\rm(iii)$ in Theorem \ref{gen2},
then the interior of $\Sigma$ is locally $G$-invariant.
\end{theorem}

\begin{remark}
Similar to Theorem \ref{localthm}, we also have local versions of Theorem \ref{immersion} and case (iv) in Theorem \ref{gen2}. Here for simplicity we omit them.
\end{remark}

\begin{proof}[Proof of Theorem \ref{localthm}]
Consider any one-parameter subgroup in $G$ and denote the infinitesimal generator by $\phi\in \mathcal{G}$. Let $\nu$ be a unit normal vector field of $\Sigma$ in $\mathbb{R}^{n+1}$ and $\Phi(X) := \langle \phi X, \nu \rangle$ be the normal component of $\phi X$ on $\Sigma$. In the following, we will  show that $\Phi(X)\equiv0$ in $\Sigma$.

In the first two cases, the mean curvature $H$ is constant. Since $\phi X$ is a Killing vector field,
according to the second variational formula of the functional $J(\Sigma)=A(\Sigma)+\frac{n}{n+1} H_0\int_{\Sigma}\langle X, \nu\rangle d\sigma$ we know that
$\Phi$ satisfies
\begin{equation}
L[\Phi]: =\Delta_{\Sigma} \Phi+|A|^2\Phi=0,\label{C1MCEL}
\end{equation}
where $\Delta_{\Sigma}$ is the Laplacian on $\Sigma$ with respect to the induced metric and $A$ is the second fundamental form of $\Sigma$ with respect to $\nu$.

As for the third case, first recall that a hypersurface in $\mathbb{R}^{n+1}$ has constant $r$-th mean curvature  if and only if it is a critical point of the functional $J_{r-1}(\Sigma):=\int_{\Sigma} (S_{r-1} +\frac{r b}{n+1}\langle X,\nu\rangle) d\sigma $ with a constant $b$.
Since $\phi X$ is a Killing vector field, according to the second variational formula of $J_{r-1}(\Sigma)$ we obtain that 
$\Phi(X)$ satisfies 
\begin{equation}
L[\Phi]:=L_{r-1}\Phi+(S_1 S_r-(r+1)S_{r+1})\Phi=0, \label{CrMCEL}
\end{equation}
where $L_{r-1}$ is defined by \eqref{LrDef} and $S_r$ is the $r$-th elementary symmetric polynomial of the principal curvatures.

It is easy to see that when $r=1$ \eqref{CrMCEL} coincides with \eqref{C1MCEL}, so there is no confusion about the operator $L$.

For each point $p\in U$, let $\{T_i\}$ be a local orthonormal frame of $G\cdot p$ around $p$ and extend it to a local orthonormal frame $\{T_i,S_j\}$ of $U$ around $p$. Then we have
\begin{equation}\label{eq5.2}
\begin{aligned}
\partial_{\mathfrak{n}}\Phi|_{p}&=\langle\phi \mathfrak{n},\nu\rangle+\langle\phi X, \bar\nabla_{\mathfrak{n}} \nu\rangle \\
=&\langle\phi(-\sin\theta \frac{X}{\|X\|} +\cos \theta N),\cos \theta \frac{X}{\|X\|} +\sin\theta N\rangle+\langle\phi X, \bar\nabla_{\mathfrak{n}} \nu\rangle \\
=&\sum_{i}\langle\phi X,T_i\rangle\langle T_i, \bar\nabla_{\mathfrak{n}} \nu\rangle+\sum_{j}\langle\phi X,S_j\rangle\langle S_j, \bar\nabla_{\mathfrak{n}} \nu\rangle\\
& +\langle\phi X,\nu\rangle\langle\nu, \bar\nabla_{\mathfrak{n}} \nu\rangle+\langle\phi X,\mathfrak{n}\rangle\langle\mathfrak{n}, \bar\nabla_{\mathfrak{n}} \nu\rangle\\
=&\sum_{i}\langle\phi X,T_i\rangle\langle T_i, \bar\nabla_{\mathfrak{n}} \nu\rangle,
\end{aligned}
\end{equation}
where $\langle\,,\,\rangle$ denotes the standard Euclidean inner product and $\bar\nabla$ denotes the Levi-Civita connection of $\mathbb{R}^{n+1}$. Here we use the fact that  $\phi X$ is tangent to the orbit $G\cdot p$ and $\phi\in\mathcal{G}\subset\mathfrak{so}(n+1)$ is anti-symmetric.

By considering the parallel hypersurfaces along the normal exponential map of $U$ in $\Sigma$ with respect to the inward unit conormal $-\mathfrak{n}$, we can extend $\{T_i, S_j, \mathfrak{n}\}$ to a local orthonormal frame of $\Sigma$ on a neighborhood $\tilde{U}$ of $p$. Moreover, at each point $q\in \tilde{U}$, $\{T_i, S_j\}$ are tangent to the parallel hypersurface of $U$ through $q$ and $\mathfrak{n}$ is normal to this parallel hypersurface.  Let $[\,,\,]$ denote the Lie bracket of vector fields along $\Sigma$. By definition $[T_i,\mathfrak n]|_p\in T_p\Sigma$ at each point $p\in U$, hence we have $\langle\nu,[T_i,\mathfrak n]\rangle|_{U}=0$. Moreover, since the contact angle $\theta$ is locally $G$-invariant, we have $T_i\theta|_U =0$. Now it follows that 
\begin{equation}\label{eq5.3}
\begin{aligned}
\langle T_i,\bar\nabla_{\mathfrak{n}} \nu\rangle&=-\langle\nu,\bar\nabla_{\mathfrak{n}}T_i\rangle\\
&=-\langle\nu,\bar\nabla_{T_i}\mathfrak{n}\rangle+\langle\nu,[T_i,\mathfrak{n}]\rangle\\
&=-\langle \cos \theta \frac{X}{\|X\|} +\sin\theta N,\bar\nabla_{T_i}(-\sin\theta \frac{X}{\|X\|} +\cos \theta N)\rangle\\
&=\sin\theta \cos \theta\langle \frac{X}{\|X\|} ,\bar\nabla_{T_i}\frac{X}{\|X\|} \rangle+\sin^2\theta\langle N,\bar\nabla_{T_i}\frac{X}{\|X\|} \rangle\\
&\quad-\cos ^2\theta\langle \frac{X}{\|X\|} ,\bar\nabla_{T_i}N\rangle-\sin\theta \cos \theta \langle N, \bar\nabla_{T_i} N\rangle\\
&=\langle N, \bar\nabla_{T_i}\frac{X}{\|X\|}\rangle\\
&=0.
\end{aligned}
\end{equation}

By Theorem \ref{Analyticity}, $X(\Sigma)$ is real analytic at each interior point and it can be extended analytically across $U$. Now taking \eqref{CrMCEL} , \eqref{eq5.2} and \eqref{eq5.3} into account, we deduce that $f=\Phi$ is a real analytic solution to the following Cauchy problem
\begin{equation}\label{Cauchy2}
\left\{\begin{array}{ll}
L[f]=0 &\textrm{in}\, \Sigma,\\
f=\partial_{\mathfrak{n}}f=0& \textrm{on}\, U.
\end{array}\right.
\end{equation}
By the Cauchy-Kovalevskaya theorem, the Cauchy problem \eqref{Cauchy2} has a unique real analytic solution (cf. Appendix A). Since there is a trivial solution $f\equiv0$, we have $f=\Phi\equiv0$ on $\Sigma$. Therefore the interior of $\Sigma$ is invariant under the infinitesimal action $\phi$, that is, $\Sigma$ is locally $G$-invariant.
\end{proof}

With the help of the preparation above, we are in position to prove Theorem \ref{gen2}.

\begin{proof}[Proof of Theorem \ref{gen2}, case (i)-(iii)]
By assumption, $\Sigma$ is complete with respect to the induced metric and each connected component of $\partial\Sigma$ is a $G$-invariant submanifold in $\mathbb{R}^{n+1}$. Moreover, we have proved in Theorem \ref{localthm} that $\Sigma$ is locally $G$-invariant. Hence, it follows from Proposition \ref{Killing} that $\Sigma$ is $G$-invariant.
\end{proof}
\begin{remark}
In the case of hypersurfaces of constant $r$-th mean curvature, we assume an additional condition that there exists an interior elliptic point in $X(\Sigma)$, that is, a point where all principal curvatures have the same sign. In particular, when $\Sigma$ is compact and $\partial\Sigma$ is a round $(n-1)$-sphere, this condition holds naturally except that $\Sigma$ is a flat disk. However, in general the compactness cannot ensure the existence of an interior elliptic point.
\end{remark}

\section{Symmetry of Helfrich-type hypersurface}\label{sec6}
We now turn to 
interior symmetry of a new class of hypersurfaces, called Helfrich-type hypersurfaces, which can be regarded as an extension of Willmore surfaces.
For an immersion $X:\Sigma^n\to\mathbb{R}^{n+1}$, the \emph{Helfrich energy} is defined by
\begin{equation*}
\mathcal{H}_{c}[\Sigma]:=\int_{\Sigma}(H+c)^2\textrm{d}\Sigma,
\end{equation*}
where $c\in\mathbb{R}$ is a constant. Obviously, when $n=2$ and $c=0$, $\mathcal{H}_0$ is exactly the Willmore energy in $\mathbb{R}^3$, and critical surfaces are just Willmore surfaces. When $n=2$, as critical surfaces of Helfrich energy, Helfrich surfaces have also been widely studied (\cite{DDG21, Helfrich73, PP21, PP22, TO04}, etc. and references therein).

Firstly, we derive the following first variational formula of the Helfrich energy $\mathcal{H}_c$ for $n\ge 2$, which contains the $n=2$ case in \cite{TO04}.

\begin{proposition}
Let $X:\Sigma^n\to\mathbb{R}^{n+1}$ be an immersion of a hypersurface with boundary $\partial\Sigma$. Consider any sufficiently smooth variation of $X$ with compactly supported variational vector field $V+\Phi \nu$, where $V$ is tangent to $X(\Sigma)$ and $\nu$ is the unit normal vector field of $X(\Sigma)$. Then the first variation of the Helfrich energy $\mathcal{H}_c$ is given by
\begin{equation*}\label{VarHelfrich}
    \begin{aligned}
    \frac{\partial}{\partial t} \mathcal{H}_c=&-\frac{2}{n}\int_{\Sigma}\{\Delta_{\Sigma}H+(H+c)[|A|^2-\frac{n^2}{2}(H+c) H]\}\Phi \mathrm{d}\Sigma\\
    &+\int_{\partial\Sigma}\{(H+c)^2\langle V,\mathfrak{n}\rangle+\frac{2}{n}[\Phi\partial_{\mathfrak{n}} H-(H+c)\partial_{\mathfrak{n}} \Phi]\} \mathrm{d}s,
    \end{aligned}
\end{equation*}
where $\mathfrak{n}$ is the outward unit conormal of $\partial \Sigma$ in $\Sigma$.
\end{proposition}
\begin{proof}
Let $X_t:=X+t(V+\Phi \nu)$. Under a natural coordinate system $\{x_i\}_{1\le i \le n}$, we have $g_{ij}=\langle \frac{\partial}{\partial x_i} X_t, \frac{\partial}{\partial x_j}X_t\rangle$ and $
    h_{ij}=-\langle \frac{\partial}{\partial x_j}\frac{\partial}{\partial x_i}X_t,\nu\rangle$. We denote the Levi-Civita connection on $\Sigma$ by $\nabla$.
    
Firstly, let us consider the case $V\equiv 0$. By direct calculation we have the following evolution equations:
\begin{equation}\label{evoEq}
    \begin{aligned}
    \frac{\partial}{\partial t} X_t &=\Phi \nu,\\
    \frac{\partial}{\partial t} g_{ij}&=2\Phi h_{ij},\\
    \frac{\partial}{\partial t} \mathrm{d}\Sigma&=n\Phi H\mathrm{d}\Sigma,\\
    \frac{\partial}{\partial t} h_{ij}&=-\nabla_i\nabla_j\Phi+\Phi h_{ik}h_{jl}g^{kl},\\
    \frac{\partial}{\partial t} H &= -\frac{1}{n}(\Delta_\Sigma \Phi+|A|^2\Phi),
    \end{aligned}
\end{equation}
 at $t=0$.
Then the first variation of $\mathcal{H}_c$ is
\begin{equation}\label{H1}
    \begin{aligned}
    \frac{\partial}{\partial t} \mathcal{H}_c=&\int_\Sigma 2(H+c)\frac{\partial}{\partial t} H\mathrm{d}\Sigma +(H+c)^2 \frac{\partial}{\partial t}\mathrm{d}\Sigma\\
    =&\int_\Sigma \{-\frac{2}{n}(H+c)(\Delta_\Sigma\Phi+|A|^2\Phi)+n(H+c)^2H\Phi \}\mathrm{d}\Sigma\\
    =&-\frac{2}{n}\int_\Sigma \{\Delta_\Sigma(H+c)+(H+c)|A|^2-\frac{n^2}{2}(H+c)^2H\}\Phi \mathrm{d}\Sigma\\
    &-\frac{2}{n}\int_{\partial\Sigma}[(H+c)\partial_{\mathfrak{n}}\Phi-\Phi\partial_{\mathfrak{n}}(H+c)]\mathrm{d}s,
    \end{aligned}
\end{equation}
where we use the Green's second identity in the last equality.

Secondly, for the case $\Phi=0$, by direct calculation we have
\begin{equation*}
    \begin{aligned}
    \frac{\partial}{\partial t} H=&\langle \nabla H,V\rangle,\\
    \frac{\partial}{\partial t} \mathrm{d}\Sigma=&\mathrm{div} V \mathrm{d}\Sigma,\\
    \end{aligned}
\end{equation*}
at $t=0$. 
Then the first variation of $\mathcal{H}_c$ is
\begin{equation}\label{H2}
    \begin{aligned}
    \frac{\partial}{\partial t} \mathcal{H}_c 
    =&\int_\Sigma [2(H+c)\langle \nabla H,V\rangle+(H+c)^2 \mathrm{div} V]\mathrm{d}\Sigma\\
    =&\int_\Sigma\mathrm{div} [(H+c)^2 V]\mathrm{d}\Sigma\\
    =&\int_{\partial\Sigma}(H+c)^2\langle V,\mathfrak{n}\rangle \mathrm{d}s.
    \end{aligned}
\end{equation}
Finally, the conclusion follows from \eqref{H1} and \eqref{H2}.
\end{proof}

\begin{definition}\label{defHelfrich}
Let $X:\Sigma^n\to\mathbb{R}^{n+1}$ be an immersion of a hypersurface with boundary $\partial\Sigma$. We say $\Sigma$ is \emph{Helfrich-type} if
\begin{equation}\label{EL}
\begin{aligned}
\Delta_{\Sigma}H+(H+c)[|A|^2-\frac{n^2}{2}(H+c) H]=0 \textrm{ in } \Sigma,\\
H=-c,\quad\partial_{\mathfrak n}H=0 \textrm{ on } \partial\Sigma.
\end{aligned}
\end{equation}
\end{definition}

Similar to Theorem \ref{Analyticity}, we have the following regularity result for Helfrich-type hypersurfaces.
\begin{theorem}\label{Analyticity2}
Let $X:\Sigma^n\to\mathbb{R}^{n+1}$  be a $C^{4,\alpha}$ embedding of a hypersurface
with boundary $\partial\Sigma$ for some $\alpha \in (0,1)$. If $\Sigma$ is a Helfrich-type hypersurface, then $\Sigma$ is real analytic at each interior point and can be extended analytically across a real analytic portion of the boundary $\partial\Sigma$.
\end{theorem}
\begin{proof}
As in the proof of Theorem \ref{Analyticity}, we write $X$ as a graph of $u$ locally. Then the Euler-Lagrange equation \eqref{EL} can be written as an equation of $u$:
\begin{equation}\label{HEL}
\begin{aligned}
\Delta_{\Sigma}\Delta_{\Sigma}u&+\frac{1}{W}(\Delta_{\Sigma}W \Delta_{\Sigma}u+\langle\nabla W, \nabla \Delta_{\Sigma}u\rangle)\\
&+ (c+\frac{W}{n} \Delta_{\Sigma}u)[\frac{n}{W^3}|\nabla^2 u|^2-\frac{n^2}{2} \Delta_{\Sigma}u (c +\frac{W}{n}\Delta_{\Sigma}u)]=0.
\end{aligned}
\end{equation}
Notice that this quasilinear equation  is of fourth order. Thus the principal part of its linearization is just the terms of order four and its characteristic polynomial becomes
\begin{equation*}
L^0(x,\lambda)=\frac{1}{W^4}\sum_{i,j,k,l=1}^n(\delta_{kl}-\frac{u_k u_l}{W^2})(\delta_{ij}-\frac{u_i u_j}{W^2})\lambda_i\lambda_j\lambda_k\lambda_l,
\end{equation*}
 where $\lambda=(\lambda_1,\cdots,\lambda_n)\in \mathbb{R}^{n}$.
Hence, the equation \eqref{HEL} is a quasilinear real analytic elliptic equation and the conclusion follows from Theorem \ref{Morrey1} and Theorem \ref{Morrey2} as before.
\end{proof}

Now based on the Euler-Lagrange equation \eqref{EL} and \eqref{HEL}, we adopt 
 the method in the proof of the first three cases in Theorem \ref{gen2} to prove the last case.

\begin{proof}[Proof of Theorem \ref{gen2}, Case (iv)]
Let $\mathcal{G}\subset \mathfrak{so}(n+1)$ be the Lie algebra of $G$.
Consider any action in $G$ and denote the infinitesimal generator by $\phi\in\mathcal{G}$. Let $\nu$ be a unit normal vector field of $\Sigma$ in $\mathbb{R}^{n+1}$ and $\Phi(X) := \langle \phi X, \nu \rangle$ be the normal component of $\phi X$. We will show that $\Phi(X)\equiv0$ in $\Sigma$.

By Theorem \ref{Analyticity2}, we see that $X(\Sigma)$ is real analytic at each interior point and it can be extended analytically across $\Gamma$. Moreover,  the first equation in \eqref{EL} can be regarded as a real analytic elliptic equation $F[X]=0$. Since $\phi X$ is a Killing vector field, $\Phi$ satisfies the linearization of this equation:
\begin{equation*}
L[\Phi]:=\frac{\partial}{\partial t}F[X+t\Phi \nu]|_{t=0}=\frac{\partial}{\partial t}F[X+t\phi X]|_{t=0}=0,
\end{equation*}
at a solution of $F[X]=0$.
In order to compute $L[\Phi]$, let us consider the normal variation $X_t:=X+t\Phi\nu$.
By using the evolution equations in \eqref{evoEq}, we can further compute the evolution of $|A|^2$
\begin{equation}\label{evoEq2}
    \frac{\partial}{\partial t} |A|^2=-2g^{ik}g^{jl}h_{kl}\nabla_i\nabla_j\Phi-2\Phi h_{ij}h_{kl}h_{mn}g^{il}g^{jm}g^{kn}.
\end{equation}
Moreover, we have
\begin{equation}\label{defLinear}
    L[\Phi]=\frac{\partial}{\partial t} F[X_t]|_{t=0}= \frac{\partial}{\partial t}\{\Delta_{\Sigma}H+(H+c)[|A|^2-\frac{n^2}{2}(H+c) H]\}|_{t=0}.
\end{equation} 
By substituting \eqref{evoEq} and \eqref{evoEq2} into \eqref{defLinear}, we  see that $L[\Phi]=0$ is a quasilinear equation of fourth order and the principal part of its linearization is just $-\frac{1}{n}\Delta_\Sigma\Delta_\Sigma \Phi$. Obviously $L[\Phi]$ is elliptic with real analytic coefficients.

For each point $p\in\Gamma$, let $\{T_i\}$ be a local orthonormal frame of  $\Gamma$ around $p$. As in the proof of Theorem \ref{gen2}, we have
\begin{equation}\label{eq:partial-Phi}
\Phi|_{\Gamma}=\partial_{\mathfrak{n}}  \Phi|_{\Gamma}=\langle\phi X, \mathfrak{n}\rangle|_{\Gamma}\equiv0,
\end{equation}
which also means $\nabla \Phi|_{\Gamma}\equiv0$.

Recall that the envolution  equation of $H$ under the variation $X_t=X+t\phi X$ reads 
\begin{equation}\label{eq4.8}
\frac{\partial}{\partial t} H= -\frac{1}{n}(\Delta_\Sigma\Phi+|A|^2 \Phi)+\langle \nabla H, \phi X\rangle.
\end{equation}
Due to the facts that $H$ is constant on $\Gamma$ and $\phi X$ is tangent to $\Gamma$, \eqref{eq4.8} 
 yields to
\begin{equation*}
0=\Delta_\Sigma\Phi|_{\Gamma}=\langle\nabla_{\mathfrak{n}}\nabla \Phi, {\mathfrak{n}}\rangle|_{\Gamma}.
\end{equation*}
Now we have
\begin{equation}\label{eq:partial2-Phi}
\begin{aligned}
\partial_{\mathfrak{n}}^2 \Phi|_{\Gamma}=\langle \nabla \Phi, \nabla_{\mathfrak{n}}\mathfrak{n}\rangle|_{\Gamma}+\langle\nabla_{\mathfrak{n}}\nabla \Phi, {\mathfrak{n}}\rangle|_{\Gamma}=0,
\end{aligned}
\end{equation}
which also means $\nabla\nabla\Phi|_{\Gamma}\equiv0$.

Since $H$ is constant and $\partial_{\mathfrak{n}} H=0$ on $\Gamma$, it follows that $\nabla H =\partial_{\mathfrak{n}}H \mathfrak{n} =0$ on $\Gamma$. 
By \eqref{eq4.8} again, on $\Gamma$ we have
\begin{equation*}
\begin{aligned}
0&=\partial_{\mathfrak{n}}(\Delta_\Sigma\Phi)+\langle\nabla_\mathfrak{n}\nabla H, \phi X\rangle\\
&=\partial_{\mathfrak{n}}(\Delta_\Sigma\Phi)+\sum_{i=1}^{n-1}\langle\nabla_\mathfrak{n}\nabla H,T_i\rangle\langle T_i, \phi X\rangle\\
&=\partial_{\mathfrak{n}}(\Delta_\Sigma\Phi)+\sum_{i=1}^{n-1}\langle\nabla_{T_i}\nabla H,\mathfrak{n}\rangle\langle T_i, \phi X\rangle\\
&=\partial_{\mathfrak{n}}(\Delta_\Sigma\Phi).
\end{aligned}
\end{equation*}
On the other hand,
\begin{equation*}
\partial_{\mathfrak{n}}(\Delta_\Sigma\Phi)=\partial_{\mathfrak{n}}^3\Phi-\nabla_{\mathfrak{n}}\langle\nabla\Phi, \nabla_{\mathfrak{n}}\mathfrak{n}\rangle+\sum_{i=1}^{n-1}\langle\nabla_{\mathfrak{n}}\nabla_{T_i}\nabla\Phi, T_i\rangle+\sum_{i=1}^{n-1}\langle\nabla_{T_i}\nabla \Phi, \nabla_{\mathfrak{n}} T_i\rangle.
\end{equation*}
It follows that
\begin{equation}
\begin{aligned}
\partial_{\mathfrak{n}}^3\Phi|_{\Gamma}&=-\sum_{i=1}^{n-1}\langle\nabla_{\mathfrak{n}}\nabla_{T_i}\nabla\Phi, T_i\rangle|_{\Gamma}\\
&=\sum_{i=1}^{n-1}\langle R(\mathfrak{n}, T_i)\nabla\Phi-\nabla_{T_i}\nabla_{\mathfrak{n}}\nabla\Phi-\nabla_{[\mathfrak{n},T_i]}\nabla\Phi,T_i\rangle|_{\Gamma}\\
&=0,
\end{aligned}
\end{equation}
where $R$ is the curvature tensor on $\Sigma$.

Now we have a real analytic solution $f=\Phi$ to the following Cauchy problem
\begin{equation}\label{Cauchy3}
\left\{\begin{array}{ll}
L[f]=0, &\textrm{in}\, \Sigma,\\
f=\partial_{\mathfrak{n}}f=\partial_{\mathfrak{n}}^2f=\partial_{\mathfrak{n}}^3 f=0, & \textrm{on}\, \Gamma.
\end{array}\right.
\end{equation}
By the Cauchy-Kovalevskaya theorem, the Cauchy problem \eqref{Cauchy3} has a unique real analytic solution (cf. Appendix A). Since there is a trivial solution $f\equiv0$, we have $f=\Phi\equiv0$ on $\Sigma$. Therefore the interior of $\Sigma$ is invariant under the infinitesimal action $\phi$, that is, $\Sigma$ is locally $G$-invariant.

Moreover, if $\Sigma$ is complete with respect to the induced metric and each connected component of $\partial\Sigma$ is a $G$-invariant submanifold in $\mathbb{R}^{n+1}$, then it follows from Proposition \ref{Killing} that $\Sigma$ is $G$-invariant.
\end{proof}
\begin{remark}
\begin{itemize}
    \item[(i)] In particular, the assumption that \emph{$\Gamma$ is real analytic and $G$-invariant} can also be replaced by that \emph{$\Gamma$ is an orbit of $G$}.
    \item[(ii)] When $n=2$, the rotational symmetry of Helfrich surface was studied by Palmer and P\'ampano in \cite{PP22}.
\end{itemize}
\end{remark}

\section{Further discussion on immersions}\label{sec7}
In this section we consider immersions of hypersurfaces. Using the fact that an immersion is locally an embedding, we can apply the proof of Theorem \ref{localthm} to obtain Theorem \ref{immersion}.

\begin{proof}[Proof of Theorem \ref{immersion}]
Consider any one-parameter subgroup in $G$ and denote the infinitesimal generator by $\phi\in\mathcal{G}$. For each point $p$ in $\Sigma$, locally we regard $X$ as an embedding in a neighborhood $U_p$. Let $\nu$ be a unit normal vector field and $\Phi(X) := \langle \phi X, \nu \rangle$ be the normal part of $\phi X$ defined on $X(U_p)$. Firstly, we will show that $\Phi\equiv0$ in any such neighborhood $X(U_p)$.

For any $p\in\partial\Sigma$, there exists a neighborhood $U_p$ such that $X|_{U_p}:U_p\to\mathbb{R}^{n+1}$ is an embedding. As in the proof of Theorem \ref{gen2}, we have $\Phi\equiv0$ in $X(U_p)$. Consider a subset $\Omega$ of $\Sigma$ defined by
\begin{equation*}
\begin{aligned}
\Omega:=\{&p\in\Sigma: \textrm{ there exists a neighborhood } U_p \\
&\textrm{ such that } X|_{U_p} \textrm{ is an embedding, and } \Phi\equiv0 \textrm{ in } X(U_p)\}.
\end{aligned}
\end{equation*}
Then $\partial\Sigma\subset \Omega$ and $\Omega$ is obviously open. For each point $q\in\partial \Omega$, let $U_q$ be a neighborhood of $q$ such that $X|_{U_q}$ is an embedding. By the definition of $\partial \Omega$, there exists a sequence of points $\{q_i\}\in \Omega \cap U_q$ such that $\lim_{n\to\infty}q_i=q$. Since $\Phi$ is real analytic in $X(U_q)$ and $\Phi (X(q_i))=0$ for each $i$, we have $\Phi\equiv 0$ in $X(U_q)$ which means $q\in \Omega$. Therefore, $\Omega$ is closed and thus $\Omega=\Sigma$. 

Now for each $p\in\Sigma$, there exists a neighborhood $U_p$ such that $X|_{U_p}$ is an embedding and $\Phi\equiv0$ in $X(U_p)$. Then we have a partial action 
$(\psi_p, \mathcal{D}_p)$ defined by
\begin{equation*}
    \begin{aligned}
    \psi_p: \mathcal{D}_p&\to X(U_p)\\
    (g,X(s))&\to \psi_p(g,X(s)):= g\cdot X(s)
    \end{aligned}
\end{equation*}
in a neighborhood $\mathcal{D}_p\subset G\times X(U_p)$ of $\{e\}\times X(U_p)$. If $X(U_p)\cap X(U_q)$ is not empty, then for each $(g,X(s))$ in $\mathcal{D}_p$ we get $\psi_p(g,X(s))=g\cdot X(s)=\psi_q(g,X(s))$. So we have a partial action $(\psi,\mathcal{D})$ in $\mathcal{D}=\bigcup_p \mathcal{D}_p$ such that $\psi|_{\mathcal{D}_p}=\psi_p$ and the interior of $X(\Sigma)$ is locally $G$-invariant.

Moreover, if $X(\Sigma)$ is closed and $X(\partial\Sigma)$ is $G$-invariant, then $X(\Sigma)$ is actually $G$-invariant. To see this, for each $p\in\mathrm{Int}\Sigma$, we only need to show that $G\cdot X(p)\subset X(\Sigma)$.  For each $\phi\in\mathcal{G}$, we denote $\theta_{\phi}(t,p)$ to be the global flow generated by $\phi X$ in $\mathbb{R}^{n+1}$. Since $G$ is a compact connected Lie group, we have 
$$
G\cdot X(p)=\bigcup_{\substack{\phi\in\mathcal{G} \\ t\in\mathbb{R}}}\theta_{\phi}(t,p).
$$

Now let $W:=\{t\in\mathbb{R}|\theta_{\phi}(t,p)\in X(\Sigma)\}$. For each $t_0\in W$, we can find $q\in\Sigma$ such that $X(q)=\theta_{\phi}(t_0,p)\in X(\Sigma$). Then there exists a neighborhood $U_q$ of $q$ such that $X|_{U_q}$ is an embedding.  
By the argument above, 
we have $\Phi\equiv 0$ in $X(U_q)$. Therefore, $\phi X$ is a Killing vector field on $X(U_q)$ and there exists an integral curve $\theta_{X(q)}(t):(-\epsilon, \epsilon)\to X(U_q)$. By the uniqueness of integral curve, we have $(t_0-\epsilon,t_0+\epsilon)\subset W$,  which means $W$ is open.

On the other hand, let $\{t_i\}$ be a sequence in $W$ converging to $b$, then $\theta_{\phi}(t_i,p)$ is a sequence in $X(\Sigma)$ converging to $\theta_{\phi}(b,p)$. Since $X(\Sigma)$ is closed in $\mathbb{R}^{n+1}$, $\theta_{\phi}(b,p)$ is also contained in $X(\Sigma)$, that is, $b\in W$. Hence $W$ is closed.
Consequently we have $W=\mathbb{R}$ and $G\cdot X(\Sigma)\subset X(\Sigma)$.
\end{proof}

\appendix
\section{The Cauchy-Kovalevskaya theorem}\label{secA}
To be complete, in this section we include the classical Cauchy-Kovalevskaya theorem,  mainly following  Rauch's book \cite{GTM128}.

In $\mathbb{R}^{n+1}$ with coordinate $(t,x^1,\cdots,x^n)$, consider the fully nonlinear partial differential equation
\begin{equation}\label{CKeq}
F(t,x,\partial^m_tu,\partial^j_t\partial^{\alpha}_xu\,|\,j\le m-1,\,  j+|\alpha|\le m)=0
\end{equation}
with prescribed data
\begin{equation}\label{CKdata}
\partial^j_tu(0,x)=g_j(x),\quad0\le j\le m-1.
\end{equation}
Given $(0,\bar{x})\in\mathbb{R}^{n+1}$, suppose that there exists $\gamma\in\mathbb{R}$ such that
\begin{equation}\label{pointSol}
F(0,\bar{x},\gamma,\partial_x^{\alpha}g_j(\bar{x}))=0.
\end{equation}
In order to solve \eqref{CKeq} for $\partial^m_tu$ by the implicit function theorem, we need the following \emph{non-characteristic} condition.
\begin{definition}
We say the hypersurface $t=0$ is \emph{non-characteristic} at $(0,\bar{x})$ with respect to the solution $\gamma$ to \eqref{pointSol}, if
\begin{equation}\label{non-cha}
\frac{\partial}{\partial s}F(0,\bar{x},s,\partial_x^{\alpha}g_j(\bar{x}))\Big|_{s=\gamma}\ne0
\end{equation}
holds.
\end{definition}
More precisely, if \eqref{pointSol}-\eqref{non-cha} hold, and $F$ is real analytic near $(0,\bar{x},\gamma,\partial_x^{\alpha}g_j(\bar{x}))$, then it follows from the implicit function theorem that 
$$
\partial^m_tu(t,x)=G(t,x,\partial^j_t\partial^{\alpha}_xu\,|\,j\le m-1,j+|\alpha|\le m),
$$
for some real analytic $G$  near $(0,\bar{x},\partial^{\alpha}_xg_j(\bar{x}))$.

Now we can state the Cauchy-Kovalevskaya theorem for fully nonlinear partial differential equations.
\begin{theorem}[Cauchy-Kovalevskaya]\label{CK}
Consider the fully nonlinear partial differential equation \eqref{CKeq} with prescribed data \eqref{CKdata}. Given $(0,\bar{x})\in\mathbb{R}^{n+1}$, suppose that there exists $\gamma\in\mathbb{R}$ satisfying \eqref{pointSol}, and $F$ and $g_j$ are real analytic near $(0,\bar{x},\gamma,\partial^{\alpha}_xg_j(\bar{x}))$ and $\bar{x}$, respectively and suppose that the hypersurface $t=0$ is non-characteristic at $(0,\bar{x})$ with respect to $\gamma$. Then, there exists a  unique real analytic solution $u$ to \eqref{CKeq} realizing \eqref{CKdata} locally around $(0,\bar{x})$ and $\partial^m_tu(0,\bar{x})=\gamma$. 
\end{theorem}

The uniqueness of real analytic solution to the Cauchy problems \eqref{Cauchy2} and \eqref{Cauchy3} can be established by  the Cauchy-Kovalevskaya theorem \ref{CK}.
\begin{theorem}\label{CK2}
$\Phi\equiv0$ is the unique analytic solution to the Cauchy problems \eqref{Cauchy2} and \eqref{Cauchy3} locally. Moreover, if $\Sigma$ is connected, then $\Phi\equiv0$ in $\Sigma$.
\end{theorem}
\begin{proof}
We only prove the case in Theorem \ref{localthm}, i.e., the Cauchy problem \eqref{Cauchy2}, since the same idea applies  to other cases.

For each point $p\in U$, there is a real analytic local parametrization $(V, f)$ with $f(0)=p$. Let $(t,x^2,\cdots,x^n)$ be the coordinates of  $\mathbb{H}^n=\{(t,x^2,\cdots,x^n)\in\mathbb{R}^n|t\ge0\}$. Then locally the Cauchy problem \eqref{Cauchy2} can be written as the following Cauchy problem
\begin{equation}\label{Cauchy2R}
\left\{\begin{array}{ll}
L[u(t,x^2,\cdots,x^n)]=0 &\textrm{in}\, V,\\
u=\partial_{t}u=0& \textrm{on}\, V\bigcap\partial\mathbb{H}^n.
\end{array}\right.
\end{equation}

To apply  Theorem \ref{CK}, we only need to verify the non-characteristic condition \eqref{non-cha}.
Under the coordinates $(x^1=t,x^2,\cdots,x^n)$, the principal part of $L[u]$ can be written as $L^0[u]=\sum_{i,j=1}^n A_{ij}(x,Du)\partial_{i}\partial_{j}u$. Since  the operator $L$ is elliptic and the ellipticity is coordinate-free, the coefficient matrix $(A_{ij}(x,Du))$ is positive definite and in particular the element $A_{11}$ is positive. Hence the hypersurface $t=0$ is non-characteristic at the origin of $\mathbb{R}^n$ for  solution $u$ of \eqref{Cauchy2R}.

According to Theorem \ref{CK}, the Cauchy problem \eqref{Cauchy2R} has a unique real analytic solution locally. As $u\equiv0$ provides a trivial  solution to \eqref{Cauchy2R}, $u=\Phi\equiv 0$ in $\Sigma$. 
\end{proof}
\begin{remark}
Actually, we only use the uniqueness part of  Theorem \ref{CK}. 
Let us include a simple proof of this. 
Assume that  $u$  is a real analytic solution to \eqref{CKeq} realizing \eqref{CKdata} and $\partial^m_tu(0,\bar{x})=\gamma$. Under the non-characteristic condition \ref{non-cha}, we can utilize the implicit function theorem to get
$$
\partial^m_tu(t,x)=G(t,x,\partial^j_t\partial^{\alpha}_xu\,|\,j\le m-1,j+|\alpha|\le m),
$$
for some real analytic $G$ near $(0,\bar{x},\partial^{\alpha}_xg_j(\bar{x}))$.

From  initial data \eqref{CKdata} it follows that
$$
\partial^j_t\partial^{\beta}_xu(0,x)=\partial^{\beta}_xg_j(x), \quad0\le j\le m-1,
$$
for any fixed $\beta$.
When $k\ge m$ and $\partial^j_t\partial^{\beta}_xu(0,x)$ is known for all $j\le k-1$, one has
\begin{equation*}
\begin{aligned}
\partial^k_t\partial^{\beta}_xu(t,x)&=\partial^{k-m}_t\partial^{\beta}_x\partial^m_tu(t,x)\\
&=\partial^{k-m}_t\partial^{\beta}_xG(t,x,\partial^j_t\partial^{\alpha}_xu\,|\,j\le m-1,j+|\alpha|\le m).
\end{aligned}
\end{equation*}
Hence, by induction all the derivatives of $u$ at $(0,\bar{x})$ will be uniquely determined. 
\end{remark}

\begin{remark}
The non-characteristic condition always holds for elliptic operators.
\end{remark}
\begin{ack}
The authors would like to thank Prof. Yuxiang Li for his helpful discussion on the regularity theory of second-order elliptic PDE. The first and the third named authors are partially supported by NSFC (No.~11831005,~12061131014). The second named author is partially supported by NSFC (No.~11871282). The fourth named author is sponsored in part by NSFC (No.~11971352,~12022109) and wishes to express his gratitude to Tsinghua and ICTP for warm hospitality. 
\end{ack}

\end{document}